\documentclass{article}

\usepackage{amsmath, amsthm, amssymb}

\newtheorem{theorem}{Theorem}[section]
\newtheorem{lemma}[theorem]{Lemma}
\newtheorem{corollary}[theorem]{Corollary}
\theoremstyle{definition}
\newtheorem{definition}[theorem]{Definition}
\newtheorem{example}[theorem]{Example}

\theoremstyle{remark}

\numberwithin{equation}{section}

\usepackage{graphicx}
\begin{document}

\title{From Weyl-Heisenberg Frames to Infinite Quadratic Forms}

\author{Xunxiang Guo, Yuanan Diao, Xingde Dai \\
Department of Mathematics\\
University of North Carolina Charlotte\\
Charlotte, NC 28223}


\maketitle

\begin{abstract}
Let $a$, $b$ be two fixed positive constants. A function $g\in
L^2({\mathbb R})$ is called a \textit{mother Weyl-Heisenberg frame
wavelet} for $(a,b)$ if $g$ generates a frame for $L^2({\mathbb
R})$ under modulates by $b$ and translates by $a$, i.e.,
$\{e^{imbt}g(t-na\}_{m,n\in\mathbb{Z}}$ is a frame for
$L^2(\mathbb{R})$. In this paper, we establish a connection
between mother Weyl-Heisenberg frame wavelets of certain special
forms and certain strongly positive definite quadratic forms of
infinite dimension. Some examples of application in matrix algebra
are provided.
\end{abstract}

\section{Introduction}

\medskip
Let $H$ be a separable complex Hilbert space. Let $B(\mathbb{H})$
denote the algebra of all bounded linear operators on
$\mathbb{H}$. Let $\mathbb{N}$ denote the set of natural numbers,
and $\mathbb{Z}$ be the set of all integers. A collection of
elements $\{x_{j}:j\in \jmath\}$ in $\mathbb{H}$ is called a
\textit{frame} if there exist constants $A$ and $B$, $0<A\leq B<
\infty$, such that
\begin{equation}
 A\|f\|^{2}\leq \sum_{j\in \jmath}|\langle f,\
x_{j}\rangle|^{2}\leq B\|f\|^{2}
\end{equation}
for all $f\in \mathbb{H}$. The supremum of all such numbers $A$
and the infimum of all such numbers $B$ are called the
\textit{frame bounds } of the frame and are denoted by $A_{0}$ and
$B_{0}$ respectively. The frame is called a \textit{tight} frame
if $A_{0}=B_{0}$ and is called a \textit{normalized tight} frame
if $A_{0}=B_{0}=1$. Any orthonormal basis in a Hilbert space is a
normalized tight frame. However, a normalized tight frame is not
necessary an orthonormal basis. Frames can be regarded as the
generalizations of orthogonal bases of Hilbert spaces. Although
the concept of frames was introduced a long time ago
(\cite{DS,SZN}), it is only in recent years that many
mathematicians have started to study them extensively. This is
largely due to the development and study of wavelet theory and the
close connections between wavelets and frames. For a glance of the
recent development and work on frames and related topics, see
\cite{BMM,D,DB2,YM}. Among those widely studied lately are the
Weyl-Heisenberg frames (or Gabor frames). Let $a$, $b$ be two
fixed positive constants and let $T_a$ and $M_b$ be the
\textit{translation operator by $a$} and \textit{modulation
operator by $b$} respectively, i.e., $T_ag(t)=g(t-a)$ and
$M_bg(t)=e^{ibt}g(t)$ for any $g\in L^2({\mathbb R})$. For a fixed
$g\in L^2({\mathbb R})$, we say that $(g,a,b)$ generates a
\textit{Weyl-Heisenberg frame} if $\{M_{mb}T_{na}g\}_{m,n\in
{\mathbb Z}}$ forms a frame for $L^2({\mathbb R})$. We also say
that the function $g$ is a \textit{mother Weyl-Heisenberg frame
wavelet} for $(a,b)$ in this case. Furthermore, a  measurable set
$E\subset \mathbb R$ is called a \textit{Weyl-Heisenberg frame
set} for $(a,\ b)$ if the function $g=\chi_E$ generates a
Weyl-Heisenberg frame for $L^2({\mathbb R})$ under modulates by
$b$ and translates by $a$, i.e., $\{e^{i mbt}g(t-na)\}_{m,n\in
\mathbb Z}$ is a frame for $L^2({\mathbb R})$. It is known that if
$ab>2\pi$, then $(g,a,b)$ cannot generate a Weyl-Heisenberg frame
for any $g\in L^2({\mathbb R})$. On the other hand, for any $a>0$,
$b>0$ such that $ab\le 2\pi$, there always exists a function $g\in
L^2({\mathbb R})$ such that $(g,a,b)$ generates a Weyl-Heisenberg
frame \cite{D}. However, in general, for any given $a>0$, $b>0$
with $ab\le 2\pi$, characterizing the mother Weyl-Heisenberg frame
wavelets $g$ for $(a,b)$ is a difficult problem. In fact, even for
the special case, i.e., for $a=2\pi$, $b=1$ and $g=\chi_{_E}$ for
some measurable set $E$, it is still an unsolved open question. In
\cite{C}, Casazza and Kalton was able to establish an equivalence
relation between this problem and a classical problem in complex
analysis regarding the unit roots of a special kind of
polynomials. There are many other works related to this subject,
for more information please refer to
\cite{BC,CC,CCJ,CL,C12,CDH,FJ,SZ}.

\medskip
In this paper,  we will establish an equivalence relation between
mother Weyl-Heisenberg frame wavelets of certain special forms and
certain strongly positive definite quadratic forms of infinite
dimension and show some examples of application in matrix
algebras. In Section 2, we will introduce some basic concepts and
preliminary lemmas regarding frames and frame sets. In Section 3,
we will briefly discuss quadratic forms of infinite dimension. In
Section 4 we will state our main results and give the proofs. Some
examples are given in Section 5.

\section{Preliminary Lemmas on Frames}

Let $E$ be a Lebesgue measurable set with finite measure, $g$ be a
function in $L^{2}(\mathbb{R})$ and $E_{g}$ be the support of $g$.
Following the notations used in \cite{YDXDQG} and \cite{DDG}, for
any $f\in L^2(\mathbb{R})$, we will let $H^0_E(f)$ and
$H_{g}^{0}(f)$ be the following formal summations:
\begin{eqnarray}{\label{hf}}
H^0_E(f)(\xi)&=& \sum_{n \in \mathbb{Z}}\langle M_n
\frac{1}{\sqrt{2\pi}} \chi_E,f\rangle M_n\frac{1}{\sqrt{2\pi}}
\chi_E(\xi),\\
H_{g}^{0}(f)&=&\sum_{n \in \mathbb{Z}}\langle M_ng,f\rangle M_ng.
\end{eqnarray}

\medskip
Two points $x$, $y\in E$ are said to be $2\pi$-translation
equivalent if $x-y=2j\pi$ for some integer $j$. The set of all
points in $E$ that are $2\pi$-translation equivalent to a point
$x$ is called the $2\pi$-translation equivalent class of $x$ and
the number of elements in this class is denoted by $\tau(x)$.

\begin{lemma}\cite{DDG}\label{f1}
Let $E$ be a Lebesgue measurable set of positive measure such that
$\tau(x)\le M$ for all $x\in E$ for some constant $M>0$, then for
any $f\in L^2(\mathbb{R})$, $H^0_E(f)$ converges to a function in
$L^2(\mathbb{R})$ under the $L^2(\mathbb{R})$ norm topology.
Moreover,
\begin{eqnarray}
H^0_E(f)(\xi)&=&\chi_{E}(\xi)\cdot\sum_{j\in\mathbb{Z}}
(f\cdot\chi_{E})(\xi+ 2\pi j)\nonumber\\
&=& \chi_{E}(\xi)\cdot\sum_{j\in\mathbb{Z}}f(\xi+ 2\pi
j)\cdot\chi_{E}(\xi+ 2\pi j). \label{he}
\end{eqnarray}
\end{lemma}

An immediate consequence of Lemma \ref{f1} is the following
corollary.

\begin{corollary}\label{cor1}
Let $0=n_0<n_1<n_2<\cdots <n_k$ be $k+1$ fixed integers and let
$E=\cup_{j=0}^k([0,2\pi)+2\pi n_j)$. For any $f\in L^2({\mathbb
R})$, let $F$ be the $2\pi$-periodical extension of the function
$\sum_{j=0}^kf(\xi+2\pi n_j), \forall \xi\in [0,2\pi)$. Then we
have
\begin{eqnarray}
H^0_E(f)&=& F\cdot \chi_E.
\end{eqnarray}
It follows that for any $f\in L^2({\mathbb R})$, we have
\begin{equation}
\langle H_E^0f, f\rangle=\frac{1}{k+1}\|F\cdot \chi_E\|^2.
\end{equation}
\end{corollary}

\begin{proof} The first part of the corollary is straight forward and is
left to the reader to verify. To see the second part, observe
first that
$$
\|F\cdot \chi_E\|^2=(k+1)\int_0^{2\pi}\vert \sum_{j=0}^k
f(\xi+2\pi n_j)\vert^2 d\xi.
$$
On the other hand,
$$
\langle H_E^0f, f\rangle=\sum_{j=0}^k\int_{[0,2\pi)+2\pi
n_j}F(\xi)\cdot \overline{f}(\xi)d\xi.
$$
For each $j$, with a suitable substitution, we have
\begin{eqnarray*}
\int_{[0,2\pi)+2\pi n_j}F(\xi)\cdot \overline{f}(\xi)d\xi =
\int_{[0,2\pi)}F(\xi)\cdot \overline{f}(\xi+2\pi n_j)d\xi.
\end{eqnarray*}
Therefore,
\begin{eqnarray*}
\langle H_E^0f, f\rangle &=&
\int_{[0,2\pi)}F(\xi)\sum_{j=0}^k\overline{f}(\xi+2\pi
n_j)d\xi\\
&=& \int_{[0,2\pi)}F(\xi)\cdot
\overline{F}(\xi)d\xi\\
&=&\frac{1}{k+1}\|F\cdot\chi_E\|^2.
\end{eqnarray*}
\end{proof}

\begin{corollary}\label{cor2}
For fixed integers $0=n_0<n_1<n_2<\cdots <n_k$, let
$E=\cup_{j=0}^k([0,2\pi)+2\pi n_j)$, $E_n=E+2n\pi$ and
$g=\frac{1}{\sqrt{2\pi}}\chi_E$. For any $f\in L^2({\mathbb R})$,
let $F_m(\xi)$ be the $2\pi$-periodical function over $\mathbb R$
defined by $F_m(\xi)=\sum_{j=0}^kf(\xi+2\pi (n_j+m))$, $\forall
\xi\in [0,2\pi)$. We have
\begin{eqnarray}{\label{cor2e}}
\sum_{m,n\in \mathbb Z}\vert\langle
M_{m}T_{2n\pi}g,f\rangle\vert^2= \frac{1}{k+1}\sum_{n\in \mathbb
Z}\|F_n\cdot \chi_{E_n}\|^2.
\end{eqnarray}
\end{corollary}

We leave the proof of Corollary \ref{cor2} to our reader. The
following more general results are from \cite{YDXDQG} and will be
needed in proving our main result.

\begin{lemma}\label{101}
If $E_{g}$ is a Lebesgue measurable set of positive measure such
that $\tau(x)\le M$ for all $x\in E_g$ for some constant $M>0$ and
$|g|\leq b$ for some constant $b>0$ on $E_{g}$, then $H_{g}^0$
defines a bounded linear operator. Furthermore, we have
\[
H_{g}^{0}(f)=g({\xi}) \sum_{j\in\mathbb{Z}}f(\xi+2\pi
j)\overline{g}({\xi}+2\pi j).
\]
\end{lemma}

\begin{lemma}{\label{le3}}
Let $0=n_{0}<n_{1}<n_{2}<\cdots<n_{k}$ be $k+1$ integers and
$a_{0}$, $ a_{1}$, $a_{2}$, $\cdots$, $a_{k}$ be $k+1$ complex
numbers. Let $E_{j}=[0,2\pi)+2n_{j}\pi$,
$G_{n}=\bigcup_{j=0}^{k}\big(E_j+2n\pi\big)$,
$g=\sum_{j=0}^{k}a_{j}\chi_{E_{j}}$ and
$g_{n}=\sum_{j=0}^{k}a_{j}\chi_{_{(E_{j}+2n\pi)}}$. For any $f\in
L^{2}(\mathbb{R})$, let $F_{n}$ be the $2\pi$-periodical function
defined by
\begin{eqnarray*}
F_{n}(\xi)&=&\sum_{j=0}^{k}f\big(\xi+2\pi(n_{j}+n)\big)
\overline{g_{n}}\big(\xi+2\pi(n_{j}+n)\big)\\
&=&\sum_{j=0}^{k}f(\xi+2\pi(n_{j}+n))\overline{a_{j}},\ \forall
\xi \in[0,2\pi),
\end{eqnarray*}
then we have
\[
\sum_{m,n\in\mathbb{Z}}|\langle M_{m}T_{2n\pi}g,
f\rangle|^{2}=\sum_{n\in\mathbb{Z}}
\|F_{n}\cdot\chi_{[0,2\pi)}\|^{2}.
\]
\end{lemma}
\begin{proof} By Lemma \ref{101},
$$
H_{g_n}^{0}f=\sum_{m\in \mathbb{Z}}\langle M_mg_n,f\rangle M_mg_n
$$
converges to
\begin{eqnarray*}
g_n(\xi) \sum_{m\in\mathbb{Z}}f(\xi+2m\pi)\overline{g_n}
(\xi+2m\pi) = F_{n}\cdot g_n.
\end{eqnarray*}
It follows that
\begin{eqnarray*}
& &\sum_{m\in \mathbb{Z}}|\langle M_mg_n,f\rangle|^2\\
&=& \langle H_{g_n}^{0}f,f\rangle\\
&=& \int_{G_{n}}F_{n}g_n\overline{f}d\xi\\
&=&\sum_{j=0}^{k}\int_{[0,2\pi)+2\pi(n+n_{j})}F_{n}
g_n\overline{f}d\xi\\
&=&\sum_{j=0}^{k}a_{j}\int_{[0,2\pi)+2\pi(n+n_{j})}F_{n}\overline{f}d\xi\\
&=&\sum_{j=0}^{k}a_{j}\int_{[0,2\pi)}F_{n}\cdot
\overline{f}(\xi+2\pi (n+n_{j}))d\xi\\
&=&\int_{[0,2\pi)}F_{n}\cdot\sum_{j=0}^{k}a_{j}\overline{f}(\xi+2\pi(n+n_{j}))d\xi\\
&=&\int_{[0,2\pi)}F_{n}\cdot\overline{F_{n}}d\xi\\
&=&\int_{[0,2\pi)}|F_{n}|^{2}d\xi=\|F_{n}\cdot
\chi_{[0,2\pi)}\|^{2}.
\end{eqnarray*}
Therefore,
\begin{eqnarray*}
&&\sum_{m,n\in\mathbb{Z}}|\langle M_{m}T_{2n\pi}g,
f\rangle|^{2}\\
&=& \sum_{m,n\in\mathbb{Z}}|\langle M_{m}g_n,
f\rangle|^{2}\\
&=& \sum_{n\in\mathbb{Z}}\big\langle\sum_{m\in\mathbb{Z}}\langle
M_mg_n,f\rangle M_mg_n, f\big\rangle\\
&=& \sum_{n\in\mathbb{Z}}\langle H_{g_n}^{0}f,f\rangle\\
&=& \sum_{n\in\mathbb{Z}}\|F_{n}\cdot \chi_{[0,2\pi)}\|^{2}.
\end{eqnarray*}
\end{proof}

\section{A Special Quadratic Form of Infinite
Dimension}

\medskip
Let $\{x_n\}$ be a real sequence in $\ell^2(\mathbb{Z})$, i.e.,
$\{x_n\}$ is a real valued sequence such that the series
$\sum_{n\in \mathbb Z}x_n^2$ is convergent. Let
$\{a_{ij}\}_{i,j\in \mathbb Z}$ be a sequence of real numbers with
$a_{ij}=a_{ji}$ . We can formally write $A=\{a_{ij}\}$ and think
of $A$ as an infinite dimensional symmetric matrix. Similarly, we
will write $x=\{x_n\}$ and $xAx^t$ for the formal sum
$\sum_{i,j\in \mathbb Z}a_{ij}x_ix_j$. If this formal sum is
convergent for all $x=\{x_n\}\in \ell^2({\mathbb Z})$, then we
will call $xAx^t$ an \textit{infinite quadratic form}. Notice that
it is easy to come up with examples of $A$ such that $xAx^t$ is
not defined for some $x$ (that is, the series $\sum_{i,j\in
\mathbb Z}a_{ij}x_ix_j$ is not convergent).

\medskip
\begin{definition}
We  say that an infinite quadratic form $xAx^t$ is strongly
positive definite if there exists a constant $c>0$ such that
$xAx^t\ge c \| x\|^2$ for all $x\in \ell^2({\mathbb Z})$, where
$\| x\|^2=\sum_{n\in \mathbb Z}x_n^2$. Similarly, one can define
negative definiteness.
\end{definition}

In this paper, we will be mainly dealing with a special type of
quadratic forms of infinite dimension, which we will describe
below.

\medskip
Let $0=n_{0}<n_{1}<n_{2}<\cdots<n_{k}$ be $k+1$ given integers,
$a_{0}$, $a_{1}$, $a_{2}$, $\cdots$, $a_{k}$ be $k+1$ given
nonzero real numbers and $p(z)=\sum_{j=0}^{k}a_{j}z^{n_{j}}$. Then
we can associate with $p(z)$ the following quadratic form of
infinite dimension:
\begin{equation}
\sum_{n\in
\mathbb{Z}}(a_{0}x_{n}+a_{1}x_{n+n_{1}}+a_{2}x_{n+n_{2}}+\cdots+a_{k}x_{n+n_{k}})^{2}
\end{equation}

We may write the above quadratic form formally as $x^{t}(A^tA)x$
by formally viewing $x$ as a column vector (of infinite
dimension), where $A$ is the upper triangular matrix constructed
by the $k+1$ given nonzero real numbers $\{a_{i}\}_{i=0}^{k}$. The
main diagonal of $A$ consists of $a_{0}$, the $n_{1}^{th}$
sub-diagonal consists of $a_{1}$, ..., the $n_{k}^{th}$
sub-diagonal consists of $a_{k}$ and all other entries are zeroes.
$A^{t}$ is the transpose of $A$. We leave it to our reader to
verify that for any $x=\{x_n\}\in \ell^2({\mathbb Z})$,
$x^{t}(A^tA)x$ is always convergent hence is a well-defined
infinite dimensional quadratic form. Apparently, such a quadratic
form is always semi-positive definite.

\section{Main Results and Proofs}

\begin{theorem}{\label{r6}}
Let $0=n_{0}<n_{1}<n_{2}<\cdots<n_{k}$ be $k+1$ given integers and
$a_{0}$, $a_{1}$, $a_{2}$, $\cdots$, $a_{k}$ be $k+1$ given
nonzero real numbers, then the following statements are
equivalent:

\smallskip\noindent
1. The infinite dimensional quadratic form
\[
\sum_{n\in\mathbb{Z}}(a_{0}x_{n}+a_{1}x_{n+n_{1}}+
a_{2}x_{n+n_{2}}+\cdots+a_{k}x_{n+n_{k}})^{2}
\]
is strongly positive definite;

\smallskip\noindent
2. The function $g=\sum_{j=0}^k a_j\chi_{_{F_j}}$ is a mother
Weyl-Heisenberg frame wavelet for $(2\pi,1)$, where
$F_j=[0,2\pi)+2n_j\pi$, $j=0,1,...,k$.
\end{theorem}

\begin{proof}
$2\Longrightarrow 1$: Let $x=\{x_{n}\}$ be a sequence in
$\ell^{2}(\mathbb{Z})$ and define $f(\xi)=x_{n}$ if
$\xi\in[0,2\pi)+2\pi n$, then $f\in L^{2}(\mathbb{R})$ by its
definition. Furthermore, it is easy to see that
$\|f\|^{2}=2\pi\cdot\|x\|^{2}=2\pi\sum_{m\in\mathbb{Z}}x_{m}^{2}$.
Thus, if $g=\sum_{j=1}^{k}a_{j}\chi_{_{F_{j}}}$ is a mother
Weyl-Heisenberg frame wavelet for $(2\pi,1)$, then there exists a
constant $A>0$ such that

\begin{equation}\label{ei}
\sum_{m,n\in\mathbb{Z}}|\langle M_{m}T_{2n\pi}g,f\rangle|^{2}\geq
A\|f\|^{2}=2\pi A\|x\|^2
\end{equation}
by the definition of the mother Weyl-Heisenberg frame wavelet. On
the other hand, by Lemma \ref{le3}, we also have
\[
\sum_{m,n\in\mathbb{Z}}|\langle M_{m}T_{2n\pi}g,f\rangle|^{2}=
\sum_{n\in\mathbb{Z}}\|F_{n}\cdot\chi_{[0,2\pi)}\|^{2},
\]
where $F_n$ is the $2\pi$-periodic function defined by
\begin{eqnarray}
F_{n}(\xi)=\sum_{j=0}^{k}f(\xi+2\pi(n+n_{j}))
\overline{g}(\xi+2\pi(n+n_{j})), \forall \xi\in [0,2\pi),
\end{eqnarray}
which is simply
\begin{eqnarray}
\sum_{j=0}^{k}f(\xi+2\pi(n+n_{j}))\overline{a_{j}}=
\sum_{j=0}^{k}a_{j}x_{n+n_{j}}
\end{eqnarray}
as one can easily check. It follows that
\[
\|F_{n}\cdot\chi_{[0,2\pi)}\|^{2}=2\pi\left(
\sum_{j=0}^{k}a_{j}x_{n+n_{j}}\right)^{2}
=2\pi\left(\sum_{j=0}^{k}a_{j}x_{n+n_{j}}\right)^{2},\] and
(\ref{ei}) becomes
\begin{equation}
2\pi\sum_{n\in\mathbb{Z}}\left(\sum_{j=0}^{k}
a_{j}x_{n+n_{j}}\right)^{2}\geq 2\pi A \|x\|^{2}.
\end{equation}
That is
\begin{equation}
\sum_{n\in\mathbb{Z}}\left(\sum_{j=0}^{k}a_{j}x_{n+n_{j}}\right)^{2}
\geq A \|x\|^{2}.
\end{equation}
Therefore, the infinite dimensional quadratic form
$\sum_{n\in\mathbb{Z}}\left(\sum_{j=0}^{k}a_{j}x_{n+n_{j}}\right)^{2}$
is strongly positive definite. This proves $2\Longrightarrow 1$.

\medskip
We will now prove $1\Longrightarrow 2$. That is, assuming that the
quadratic form $\sum_{n\in \mathbb
Z}\left(\sum_{j=0}^ka_{j}x_{n+n_j}\right)^2$ is strongly positive
definite, we need to prove that the function $g=\sum_{j=0}^k
a_j\chi_{_{([0,2\pi)+2n_j\pi)}}$ is a mother Weyl-Heisenberg frame
wavelet for $(2\pi,1)$. In other word, we need to prove that there
exists a constant $c_0>0$ such that
\begin{equation}\label{2to3}
\sum_{m,n\in \mathbb Z}\vert\langle
M_{m}T_{2n\pi}g,f\rangle\vert^2\ge c_0\|f\|^2
\end{equation}
for any $f\in L^2({\mathbb R})$. For any given $f\in L^2({\mathbb
R})$ and $0<\epsilon<1$, there exists $M_1>0$ such that
\begin{equation}\label{s3}
\|f\cdot \chi_{_{[-2N\pi,2M\pi]}}\|^2 \ge \epsilon\|f\|^2
\end{equation}
for any $N$, $M$ such that $N\ge M_1$ and $M\ge M_1$. Let $F_n$ be
as defined in Lemma \ref{le3}. By Lemma \ref{le3}, we have
\begin{equation}\label{s1}
\sum_{m,n\in\mathbb{Z}}|\langle M_{m}T_{2n\pi}g,
f\rangle|^{2}=\sum_{n\in\mathbb{Z}}
\|F_{n}\cdot\chi_{_{[0,2\pi)}}\|^{2}.
\end{equation}
Now choose $M=M_1+n_k$. We have
\begin{eqnarray}
& &\sum_{m,n\in \mathbb Z}\vert\langle
M_{m}T_{2n\pi}g,f\rangle\vert^2\nonumber\\
&=& \sum_{n\in\mathbb{Z}}
\|F_{n}\cdot\chi_{_{[0,2\pi)}}\|^{2}\nonumber\\
&\ge & \sum_{-M\le n\le M}
\|F_{n}\cdot\chi_{_{[0,2\pi)}}\|^{2}\nonumber\\
&=& \sum_{-M\le n\le M}\int_0^{2\pi} \big(\sum_{j=0}^k
a_jf(\xi+2\pi(n_j+n))\big)^2d\xi.\label{s2}
\end{eqnarray}
Since the quadratic form $\sum_{ n\in \mathbb{Z}}\big(\sum_{j=0}^k
a_jx_{(n+n_j)})\big)^2$ is strongly positive definite, there
exists a constant $c_1>0$ such that
\begin{equation}
\sum_{ n\in \mathbb{Z}}\big(\sum_{j=0}^k a_jx_{(n+n_j)})\big)^2
\ge c_1\|x\|^2.
\end{equation}
In particular, for each fixed $\xi\in [0,2\pi)$, we may define
$x\in \ell^2(\mathbb{Z})$ by $x_{n+n_j}= f(\xi+2\pi(n_j+n))$ if
$-M\le n\le M+n_k$ and $x_{n+n_j}=0$ otherwise. Then we have
\begin{eqnarray}
& & \sum_{-M\le n\le M}\big(\sum_{j=0}^k
a_jf(\xi+2\pi(n_j+n))\big)^2\nonumber\\
&=& \sum_{ n\in \mathbb{Z}}\big(\sum_{j=0}^k a_jx_{n+n_j})\big)^2
\ge c_1\|x\|^2\nonumber\\
&=& c_1\sum_{-M\le n\le M+n_k}|f(\xi+2\pi n)|^2.
\end{eqnarray}
It follows by (\ref{s3}) that
\begin{eqnarray*}
& &\sum_{-M\le n\le M}\int_0^{2\pi}
\big(\sum_{j=0}^k a_jf(\xi+2\pi(n_j+n))\big)^2d\xi\\
&\ge & c_1\sum_{-M\le n\le M+n_k}\int_0^{2\pi}|f(\xi+2\pi
n)|^2d\xi\\
&= & c_1\|f\cdot
\chi_{_{[-2M\pi,2(M+n_k+1)\pi)}}\|^2\\
&\ge & c_1 \epsilon\|f\|^2.
\end{eqnarray*}
\end{proof}

Notice that in the above proof, $\epsilon$ can be arbitrarily
close to $1$. It follows that the constant $c_0$ can be taken as
$c_1$. With some additional but small modifications to the above
proof, we can draw a stronger conclusion, which we state as the
following corollary. The details are left to our reader.

\begin{corollary}
Let $0=n_{0}<n_{1}<n_{2}<\cdots<n_{k}$ be $k+1$ given integers,
and $a_{0}$, $a_{1}$, $a_{2}$, $\cdots$, $a_{k}$ be $k+1$ given
nonzero real numbers, then if the infinite dimensional quadratic
form
\[
\sum_{n\in\mathbb{Z}}(a_{0}x_{n}+a_{1}x_{n+n_{1}}+
a_{2}x_{n+n_{2}}+\cdots+a_{k}x_{n+n_{k}})^{2}
\]
is strongly positive definite such that
\[
c_1\|x\|^2\le \sum_{n\in\mathbb{Z}}(a_{0}x_{n}+a_{1}x_{n+n_{1}}+
a_{2}x_{n+n_{2}}+\cdots+a_{k}x_{n+n_{k}})^{2}\le c_2\|x\|^2
\]
for some positive constants $c_1\le c_2$, then the function
$g=\sum_{j=0}^k a_j\chi_{_{F_j}}$ is a mother Weyl-Heisenberg
frame wavelet for $(2\pi,1)$ with $c_1$ as a lower frame bound and
$c_2$ as an upper frame bound. Conversely, if the function
$g=\sum_{j=0}^k a_j\chi_{_{F_j}}$ is a mother Weyl-Heisenberg
frame wavelet for $(2\pi,1)$ with $c_1$ as a lower frame bound and
$c_2$ as an upper frame bound, then
\[
c_1\|x\|^2\le \sum_{n\in\mathbb{Z}}(a_{0}x_{n}+a_{1}x_{n+n_{1}}+
a_{2}x_{n+n_{2}}+\cdots+a_{k}x_{n+n_{k}})^{2}\le c_2\|x\|^2
\]
for any $x\in \ell^2(\mathbb{Z})$.
\end{corollary}

Let us point out that by a similar approach used in \cite{C}, we
can show that the frame bounds of the Weil-Heisenberg frame
generated by the function $g$ (using $(2\pi,1)$) defined above can
be obtained by evaluating the minimum and maximum values of
$|a_0+a_1z^{n_1}+\cdots + a_kz^{n_k}|^2$ over the unit circle
$\mathbb{T}$. Let $\ell^2[M_1,M_2]$ be the subset of
$\ell^2(\mathbb{Z})$ which contain all the elements $x$ such that
$x_n=0$ for any $n<M_1$ or $n>M_2$. If we restrict $x$ to
$\ell^2[M_1,M_2]$, then the quadratic form
$\sum_{n\in\mathbb{Z}}(a_{0}x_{n}+a_{1}x_{n+n_{1}}+
a_{2}x_{n+n_{2}}+\cdots+a_{k}x_{n+n_{k}})^{2}$ is finite and its
corresponding matrix is a main diagonal block of $A$, where $A$ is
the (infinite dimensional) matrix of the quadratic form
$\sum_{n\in\mathbb{Z}}(a_{0}x_{n}+a_{1}x_{n+n_{1}}+
a_{2}x_{n+n_{2}}+\cdots+a_{k}x_{n+n_{k}})^{2}$. Conversely, any
main diagonal block of $A$ can be viewed as the matrix of a finite
quadratic form with $x$ chosen from some suitable
$\ell^2[M_1,M_2]$. This observation leads to our last theorem.

\begin{theorem}{\label{r9}}
Let $0=n_{0}<n_{1}<n_{2}<\cdots<n_{k}$ be $k+1$ given integers,
$a_{0}$, $a_{1}$, $a_{2}$, $\cdots$, $a_{k}$ be $k+1$ given
nonzero real numbers, and $A$ be the symmetrical infinite matrix
corresponding to the infinite quadratic form
$$
\sum_{n\in
\mathbb{Z}}(a_{0}x_{n}+a_{1}x_{n+n_{1}}+a_{2}
x_{n+n_{2}}+\cdots+a_{k}x_{n+n_{k}})^{2}.
$$
Let $\min_{z\in \mathbb{T}}|a_0+a_1z^{n_1}+\cdots +
a_kz^{n_k}|^2=C_{1}$ and $\max_{z\in
\mathbb{T}}|a_0+a_1z^{n_1}+\cdots + a_kz^{n_k}|^2=C_{2}$, then the
eigenvalues of any main diagonal block of $A$ are bounded between
$C_{1}$ and $C_{2}$.
\end{theorem}

\section{Examples}

A symmetric infinite dimensional matrix $B=\{b_{ij}\}$ is called
\textit{periodic} if there exists a real sequence
$\{b_k\}=\{b_0,b_1,...$ such that $b_{ij}=b_k$ whenever $|i-j|=k$.
Determining whether $B$ is strongly positive definite is a hard
question in general, since one would have to consider all the main
diagonal blocks of $B$. However, if the sequence
$\{b_k\}=\{b_0,b_1,...,b_n\}$ is a finite sequence, $B$ may be
related to a polynomial $p(z)=a_0+a_1z+\cdots+a_nz^n$ via the
following equation
\begin{eqnarray*}
\sum_0^na^2_i&=&b_0\\
\sum_0^{n-1} a_ia_{i+1}&=&b_1\\
\sum_0^{n-2} a_ia_{i+2}&=&b_2\\
...&=&...\\
a_0a_{n}&=&b_n,
\end{eqnarray*}
in which case whether $B$ is strongly positive definite can be
solved using Theorem \ref{r9}.

\begin{example}
The infinite quadratic form $\sum_{n\in
\mathbb{Z}}(x_{n}+x_{n+1}+x_{n+3})^2$ is strongly positive
definite since its corresponding polynomial $p(z)=1+z+z^3$ has the
property $min_{z\in \mathbb{T}}{|p(z)|^2}\approx 0.3689$ and
$max_{z\in \mathbb{T}}{|p(z)|^2}=9$. The symmetrical infinite
matrix corresponding to $\sum_{n\in
\mathbb{Z}}(x_{n}+x_{n+1}+x_{n+3})^2$ is
\[
A=\left(\begin{array}{ccccccccccccccc}\cdot&  &  &
\cdot&  &  & \cdot&  &  & \cdot&  &  &\cdot&  &  \\
 \cdot & 0 & 1&1&1&3&1&1&1&0&\cdot &\cdot &\cdot&\cdot&\cdot\\
\cdot & \cdot & 0 & 1&1&1&3&1&1&1&0&\cdot &\cdot &\cdot &\cdot\\
\cdot & \cdot & \cdot & 0 & 1&1&1&3&1&1&1&0&\cdot &\cdot &\cdot\\
\cdot&\cdot & \cdot & \cdot & 0 & 1&1&1&3&1&1&1&0&\cdot &\cdot \\
\cdot &\cdot &\cdot &\cdot &\cdot & 0&1&1&1&3&1&1&1&0&\cdot\\
\cdot&  &  &
\cdot&  &  & \cdot&  &  & \cdot&  &  &\cdot&  &  \\
\end{array}\right).
\]
By Theorem \ref{r9}, any eigenvalue $\lambda$ of any main diagonal
block of $A$ must satisfy
\[
0.3689\le \lambda\le 9
\]
The following are some main diagonal block matrices of the above
infinite matrix along with their eigenvalues.

\medskip
$A_2=\left(\begin{array}{cc} 3& 1
\\ 1& 3
\end{array}\right),$
eigenvalues are: $\lambda_1=4$, $\lambda_2=2$.

\medskip
$A_3=\left(\begin{array}{ccc}3& 1& 1\\1& 3& 1\\1& 1& 3
\end{array}\right),$ eigenvalues are: $\lambda_1=2$,
$\lambda_2=2$, $\lambda_3=5$.

\medskip
$A_4=\left(\begin{array}{cccc}3& 1& 1& 1\\1& 3& 1& 1\\1& 1& 3&
1\\1& 1& 1& 3
\end{array}\right),$ eigenvalues are: $\lambda_1=2$,
$\lambda_2=2$, $\lambda_3=2$, $\lambda_4=6$.

\medskip
$A_5=\left(\begin{array}{ccccc}3& 1& 1& 1& 0\\1& 3& 1& 1& 1\\1& 1&
3& 1& 1\\1& 1& 1& 3& 1\\0& 1& 1& 1& 3
\end{array}\right)$, eigenvalues are: $\lambda_1\approx 1.35$,
$\lambda_2=2$, $\lambda_3=2$, $\lambda_4=3$, $\lambda_5\approx
6.65$.
\end{example}

\begin{example}
Let $p(z)=2+3z^2+4z^3$. We have $min_{z\in
\mathbb{T}}{|p(z)|^2}=1$, $max_{z\in \mathbb{T}}{|p(z)|^2}=81$.
Thus the infinite quadratic form $\sum_{n\in
\mathbb{Z}}(2x_{n}+3x_{n+2}+4x_{n+3})^2$ is strongly positive
definite. The corresponding symmetrical infinite matrix is
\[
B=\left(\begin{array}{ccccccccccccccc}\cdot&  &  &
\cdot&  &  & \cdot&  &  & \cdot&  &  &\cdot&  &  \\
 \cdot & 0 & 8&6&12&29&12&6&8&0&\cdot &\cdot &\cdot&\cdot&\cdot\\
\cdot & \cdot & 0 & 8&6&12&29&12&6&8&0&\cdot &\cdot &\cdot &\cdot\\
\cdot & \cdot & \cdot & 0 & 8&6&12&29&12&6&8&0&\cdot &\cdot &\cdot\\
\cdot&\cdot & \cdot & \cdot & 0 & 8&6&12&29&12&6&8&0&\cdot &\cdot \\
\cdot &\cdot &\cdot &\cdot &\cdot & 0&8&6&12&29&12&6&8&0&\cdot\\
\cdot&  &  &
\cdot&  &  & \cdot&  &  & \cdot&  &  &\cdot&  &  \\
\end{array}\right).
\]
Again, any eigenvalue $\lambda$ of any main diagonal block of $B$
must satisfy
\[
1\le \lambda\le 81.
\]
A few main diagonal blocks of the above infinite matrix along with
their eigenvalues are listed below.

\medskip
$B_2=\left(\begin{array}{cc} 29& 12
\\ 12& 29
\end{array}\right),$
eigenvalues are: $\lambda_1=17$, $\lambda_2=41$.

\medskip
$B_3=\left(\begin{array}{ccc}29& 12& 6\\12& 29& 12\\6& 12& 29
\end{array}\right),$ eigenvalues are: $\lambda_1=14.77$,
$\lambda_2=23$, $\lambda_3=49.23$.

\medskip
$B_4=\left(\begin{array}{cccc}29& 12& 6& 8\\12& 29& 12& 6\\6& 12&
29& 12\\8& 6& 12& 29
\end{array}\right),$ eigenvalues are: $\lambda_1\approx 12.68$,
$\lambda_2\approx 20.89$, $\lambda_3\approx 25.32$,
$\lambda_4\approx 57.11$.

\medskip
$B_5=\left(\begin{array}{ccccc}29& 12& 6& 8& 0\\12& 29& 12& 6&
8\\6& 12& 29& 12& 6\\8& 6& 12& 29& 12\\0& 8& 6& 12& 29
\end{array}\right)$, eigenvalues are: $\lambda_1\approx 9.84$,
$\lambda_2\approx 20.69$, $\lambda_3=21$, $\lambda_4=31$,
$\lambda_5\approx 62.47$.
\end{example}

\begin{example} Consider the infinite dimensional matrix
\[
C=\left(\begin{array}{ccccccccccccccc}\cdot&  &  &
\cdot&  &  & \cdot&  &  & \cdot&  &  &\cdot&  &  \\
 \cdot & 0 & -2&1&-2&4&-2&1&-2&0&\cdot &\cdot &\cdot&\cdot&\cdot\\
\cdot & \cdot & 0 & -2&1&-2&4&-2&1&-2&0&\cdot &\cdot &\cdot &\cdot\\
\cdot & \cdot & \cdot & 0 & -2&1&-2&4&-2&1&-2&0&\cdot &\cdot &\cdot\\
\cdot&\cdot & \cdot & \cdot & 0 & -2&1&-2&4&-2&1&-2&0&\cdot &\cdot \\
\cdot &\cdot &\cdot &\cdot &\cdot & 0&-2&1&-2&4&-2&1&-2&0&\cdot\\
\cdot&  &  &
\cdot&  &  & \cdot&  &  & \cdot&  &  &\cdot&  &  \\
\end{array}\right).
\]
One can show that it is related to the polynomial $p(z)=-2+z+z^3$,
which has a unit zero. It follows that the the infinite quadratic
form $\sum_{n\in \mathbb{Z}}(-2x_{n}+x_{n+2}+x_{n+3})^2$ is not
strongly positive definite. Therefore, there must exist a sequence
of main diagonal blocks of $C$ whose least eigenvalues will
approach $0$ as the dimensions of the blocks go to infinity. The
following is the list of the first few main diagonal blocks and
their corresponding eigenvalues.

\medskip
$C_2=\left(\begin{array}{cc} 4& -2
\\ -2& 4
\end{array}\right),$
eigenvalues are: $\lambda_1=2$, $\lambda_2=6$.

\medskip
$C_3=\left(\begin{array}{ccc}4& -2& 1\\-2& 4& -2\\1& -2& 4
\end{array}\right),$ eigenvalues are: $\lambda_1\approx 1.63$,
$\lambda_2=3$, $\lambda_3\approx 7.37$.

\medskip
$C_4=\left(\begin{array}{cccc}4& -2& 1& -2\\-2& 4& -2& 1\\1& -2&
4& -2\\-2& 1& -2& 4
\end{array}\right),$ eigenvalues are: $\lambda_1=1$,
$\lambda_2=3$, $\lambda_3=3$, $\lambda_4=9$.

\medskip
$C_5=\left(\begin{array}{ccccc}4& -2& 1& -2& 0\\-2& 4& -2& 1&
-2\\1& -2& 4& -2& 1\\-2& 1& -2& 4& -2\\0& -2& 1& -2& 4
\end{array}\right)$, eigenvalues are: $\lambda_1\approx 0.22$,
$\lambda_2\approx 2.70$, $\lambda_3=3$, $\lambda_4=4$,
$\lambda_5\approx 10.08$.
\end{example}

\bibliographystyle{amsalpha}

\end{document}